\documentclass[12pt]{article}
\begin{document}
\def\id{\mathrm{id}}
\def\C{\mathbf{C}}
\def\bC{\mathbf{\overline{C}}}
\def\Sing{\mathrm{Sing}}
\def\B{\mathbf{B}}
\def\S{\mathbf{S}}

\title{Singularities of inverse functions}
\author{Alexandre Eremenko\thanks{Supported by NSF grant DMS-1067886.}}
\maketitle
\begin{abstract}
A survey of general results on the
singularities of inverses to meromorphic functions
is given, with applications to holomorphic dynamics.
This is a lecture delivered at the workshop
``The role of complex analysis in complex dynamics'' in Edinburgh
on May 22 2013, with corrected mistakes and updated references.
\end{abstract}

\noindent
{\bf 1. Definition of singularities}
\vspace{.2in}

Let $f:D\to G$ be a non-constant holomorphic map between Riemann surfaces.
Let $z_0$ be a point in $D$ such that $f'(z_0)\neq 0$. Then by the
inverse function theorem there exists a neighborhood $V$ of the
point $w_0=f(z_0)$ and a holomorphic map $\phi:V\to D$,
such that $f\circ\phi=\id_{V}$.

What happens when we perform an analytic
continuation of $\phi$? Let $\gamma:[0,1]\to G$ be a curve from $w_0$ to $w_1$
and suppose that an analytic continuation of $\phi$ 
along $\gamma$ is possible for $t\in[0,1)$, and let us see what can
happen when $t\to 1$. 

Consider the image $\Gamma(t)=\phi(\gamma(t)).$ There are two possibilities:

a) The curve $\Gamma(t)$ has a limit point $z_1\in D$ as $t\to 1$.
By continuity we have $f(z_1)=w_1$. Now we conclude that the limit
set of $\Gamma(t)$ must consist of one point, because otherwise
the limit set of the curve $\Gamma$ would contain a continuum,
while the preimage of a point
under $f$ is discrete. Thus $\Gamma$ ends at $z_1$. If $f'(z_1)\neq 0$,
then the analytic continuation of $\phi$ to $w_1$ is possible, and if
$f'(z_1)=0,$ then $\phi$ has an algebraic singularity (branch point)
at $w_1$.

b) The curve $\Gamma(t)$ tends to $\infty$, where $\infty$ is the added
point of the one-point compactification of $D$.
In this case $\Gamma$ is an {\em asymptotic curve} of $f$ which means
that $\Gamma$ is a curve in $D$ parametrized by $[0,1)$, $\Gamma(t)\to\infty$
as $t\to\infty$, and $f(\gamma(t))$ has a limit in $G$ as $t\to 1$.
The limits of $f$ along asymptotic curves are called
{\em asymptotic values}.

Thus non-algebraic singularities of the inverse function $f^{-1}$ correspond
to asymptotic curves of $f$.

To obtain a one-to-one correspondence, we have to define precisely
the notion of singularity, and to introduce some equivalence
relation on the asymptotic curves. But first we notice the following:
\vspace{.1in}

\noindent
{\bf Proposition 1.} {\em If $G$ contains no critical values and no asymptotic
values, then $f:D\to G$ is a covering map.}
\vspace{.1in}

We recall that a continuous map $f$ is called a {\em covering}
if every point $w_1$
in the image has a neighborhood $V$ such that every component of $f^{-1}(V)$
is mapped onto $V$ homeomorphically. An equivalent property is that every
path in the image has a unique lifting. In our situation this is
equivalent to saying that $\phi$ can be analytically continued along any
path in $G$.

Now we give an exact definition of a singularity of the inverse function.
Let us assume that $f(D)$ is dense in $G$.
Let us fix a point $a\in G$. Every neighborhood $V$ of $a$ has
non-empty open preimage. Consider a map $S$ which to
every neighborhood $V$ of $a$ puts into correspondence some component
$S(V)$ of $f^{-1}(V)$, so that the following condition is satisfied
$$V_1\subset V_2\longrightarrow S(V_1)\subset S(V_2).$$ 
Then there are two possibilities:

a) Intersection of all $S(V)$ is not empty. Then it must consist of one
point $z\in D$ such that $f(z)=a$. Indeed, let this intersection
contain a point $z$, and let $U$ be a neighborhood of $z$ then $f(U)$
is a neighborhood of $a$, and there exists a neighborhood $V$ of $a$
such that $\overline V\subset
f(U)$. This implies that $S(V)\subset U$. So $\cap_VS(V)$ consists
of one point $z$.

b) Intersection of all $S(V)$ is empty. Then we say that our map
$S$ defines a {\em transcendental singularity} of $f^{-1}$. We say
that the transcendental singularity $S$ {\em lies over $a$},
and that $a$ is the projection of the singularity $S$.

Introduction of transcendental singularities is a sort of completion of
$D$ to which our function $f$ extends continuously.
Let $\Sing(f)$ be the set of all transcendental singularities,
and $\overline{D}_f=D\cup\Sing(f),$ the disjoint union.
We define the topology on $\overline{D}_f$ as follows:
the neighborhoods of a point in $D$ are its usual neighborhoods,
and the neighborhoods of a point $S\in\Sing(f)$ are the sets $S(V)\cup\{ S\}$.
So the sets $S(V)$ are {\em punctured neighborhoods} of $S$.
We can define $f(S)=a$ if $S$ lies over $a$, and this extension
of $f$ to $\overline{D}_f$ is continuous.

It is easy to see that for each transcendental singularity $S$ over
$a$ there exists a curve
$\Gamma:[0,1)\to D$ such that for every neighborhood $V$ of
$a$ we have $f(\Gamma(t))\in V$ for all $t$ sufficiently close to $1$.
Thus $\Gamma$ is an asymptotic curve with asymptotic value $a$.
And conversely, if $\Gamma$ is an asymptotic curve with asymptotic
value $a$, then we choose as $S(V)$ to be that component of $f^{-1}(V)$
which contains $\Gamma(t)$ for $t$ sufficiently close to $1$,
and this defines a transcendental singularity.

If $S$ is a transcendental singularity, then all regions $S(V)$ are
unbounded. They are called {\em tracts} of $f$ over $a$.   

One can also define transcendental singularities as elements
of the completion
of $D$ with respect to some metric.
Suppose that $G$ is equipped by some
intrinsic metric $\sigma$. ``Intrinsic'' means that the distance between two
points is equal to the infimum of the lengths of curves connecting
these two points, for example, any smooth Riemannian metric is 
intrinsic. The pull-back $\rho=f^*\sigma$ is an intrinsic metric
in $D$ defined as follows:
the $\rho$-length of a curve in $D$ is the $\sigma$-length of its image.
So $f$ becomes a local isometry $(D,\rho)\to (G,\sigma)$.

Now, we define another metric $\rho_M$ in $D$,
which is called the {\em Mazurkiewicz metric}. The $\rho_M$-distance between
two points is the infimum of $\sigma$-{\em diameters} of curves
$f\circ\gamma$ over all curves $\gamma$
connecting these two points. Ma\-zur\-kie\-wicz's metric is in general
not intrinsic, and $\rho_M\leq\rho.$ However $\rho_M$ coincides
with $\rho$ on sufficiently small neighborhoods of a point $z\in D$ for
which $f'(z)\neq 0$.
\vspace{.1in}

\noindent
{\bf Example.} Let $f(z)=\cos z, \; \C\to\C.$
Let $\sigma$ be the usual Euclidean metric. Then $\rho=f^*\sigma$
is the metric whose line element is
$$|f'(z)||dz|=|\sin z||dz|.$$
The $\rho$-distance between the points $0$ and $2\pi m$ is $2m$,
(the shortest curve is the segment $[0,2\pi m]$ which is mapped by $f$
onto the segment $[-1,1]$ described $2m$ times). The Mazurkiewicz
distance between the same points is $2$.
\vspace{.1in}

\noindent
{\bf Exercise.} {\em Function $f$ extends to a continuous function on
the completion $\overline{D}_M$ of $D$ with respect to $\rho_M$,
and that there is a homeomorphism $\phi:\overline{D}_M\to\overline{D}_f$,
such that
the extension satisfies
such that $f(z)=f(\phi(z))$ for all $z\in \overline{D}_M$.}
\vspace{.2in}

\noindent
{\bf 2. Iversen's classification}
\vspace{.1in}

We call critical points and transcendental
singularities considered as points of $\overline{D_f}$ the
singularities of $f^{-1}$. 

We begin with the simplest kind of transcendental singularities,
the isolated ones.
Let $S$ be an isolated transcendental singularity over a point $a$,
then there is an open $\sigma$-disc $V=B(a,r)$ of radius $r$ around $a$,
such that $S(V)$ is at positive distance from other singularities.
Proposition 1 applied to the restriction
\begin{equation}\label{1}
f:S(V)\backslash f^{-1}(a)\to V\backslash\{ a\}
\end{equation}
implies that this restriction is a covering map.
All possible coverings over a punctured disc are classified
by subgroups of the fundamental group which is the infinite cyclic
group. Thus there are two possibilities: 

a) (\ref{1}) is $m$-to-$1$, and $S$ is a critical point, or

b) (\ref{1}) is a universal cover. In this case $S(V)$ is a simply connected
region bounded by a single curve in $D$, parametrized by $(0,1)$, and
both ends of the curve are at $\infty$. The map (\ref{1}) is
equivalent to $z\mapsto\exp(z)$ from the left half-plane to the punctured
unit disc. This type of singularity is called {\em logarithmic}.
\vspace{.1in}

\noindent
{\bf Examples.} Function $\exp:\C\to\C$ has one logarithmic
singularity over $0$.
Function $\exp:\C\to\bC$ has two logarithmic singularities,
one over $0$ another
over $\infty$. Function $\cos:\C\to\bC$ has two logarithmic singularities 
over $\infty$ and infinitely many algebraic ones over $1$ and $-1$.
The entire function $z\mapsto\sin z/z$,
$\C\to\bC$ has two logarithmic singularities
over $\infty$, infinitely many critical points, and two non-isolated
singularities over $0$.
\vspace{.1in}

Further classification of transcendental singularities is due to
Felix Iversen (1912).

A transcendental singularity $S$ over $a$
is called {\em direct} if there exists $V$
such that $f(z)\neq a$ for $z\in S(V)$. Otherwise it is
called {\em indirect}.

So logarithmic singularities are direct.
\vspace{.1in}

\noindent
{\bf Examples.} {\em We consider functions $\C\to\bC$.
Function $e^z\sin z$ has one direct non-logarithmic singularity over $\infty$.
Function $\sin z/z$ has two indirect singularities over $0$, and
function $\sin\sqrt{z}/\sqrt{z}$ has one indirect singularity over $0$.}
\vspace{.2in}

\noindent
{\bf 3. Meromorphic functions of finite order}
\vspace{.2in}

From now on we only consider maps $\C\to\bC$, and the Riemannian metric
$\sigma$
will be the spherical metric,
whose pullback $\rho=f^*\sigma$ has the length element
$$\frac{|f'(z)||dz|}{1+|f(z)|^2}.$$
The lower order of growth of a meromorphic function
will play an important role.
It is defined by
$$\lambda(f)=\liminf_{r\to\infty}\frac{\log A(r,f)}{\log r},$$
where
$$A(r,f)=\frac{1}{\pi}\int_{|z|\leq r}\frac{|f'|^2}{(1+|f|^2)^2}dm,$$
where $dm$ is the Euclidean area element in the plane. The order
of $f$ is defined similarly, with $\limsup$ instead of $\liminf$.
The geometric interpretation of the quantity $A(r,f)$ is the
``average number of sheets'' of the covering of $\bC$
by the image of the disc $|z|\leq r$: the integral is the
area of this image, and $\pi$ is the area of $\bC$.

The Nevanlinna characteristic is defined by the averaging of $A(r,f)$,
$$T(r,f)=\int_0^r A(r,f)\frac{dt}{t}.$$
It has an advantage that it satisfies the usual properties of the degree
of a rational function,
$$T(r,f+g)\leq T(r,f)+T(r,g)+O(1),\quad T(r,fg)\leq T(r,f)+T(r,g)+O(1),$$
$$T(r,f^n)=nT(r,f)+O(1)\quad\mbox{and}\quad T(r,f)=T(r,1/f).$$
The order and lower order can be defined using $A(r,f), T(r,f)$ or,
$\log M(r,f)$ in the case of entire functions, with the same result.

\vspace{.1in}

\noindent
{\bf Theorem 1.} (Denjoy, Ahlfors, Beurling, Carleman).
{\em The number of direct singularities of $f^{-1}$
is at most $\max\{1,2\lambda(f)\}$.}
\vspace{.1in}

\noindent
{\bf Corollary.} {\em If $f$ is an entire function, then
the number of transcendental
singularities over points in $\C$ is at most $2\lambda$.
In particular, the number of finite
asymptotic values of an entire function is at most $2\lambda(f)$.}
\vspace{.1in}

To derive the Corollary, one notices that the number of singularities
of an entire function over infinity is at least the number of transcendental
singularities over finite points. Indeed, between any two asymptotic
curves corresponding to distinct transcendental singularities
there must be an asymptotic
curve with infinite asymptotic value. All singularities
over infinity are direct.

There are several different proofs of Theorem~{1}, 
due to Ahlfors, Beurling and Carleman; each of them introduced new 
important tools of analysis.
\vspace{.1in}

The main result used in Carleman's proof of Theorem 1
is the following important inequality
from potential theory.
Let $u$ be a non-negative subharmonic function
in a ring, $A=\{ z:r_0<|z|<r_1\}.$
Let
$$m_2(r)=\left(\int_0^{2\pi}u^2(re^{i\phi})d\phi\right)^{1/2},$$
$$\mu(t)=m_2(e^t).$$
Then for $r\in(r_0,r_1)$ we have
\begin{equation}\label{carleman}
\mu^{\prime\prime}(t)\geq \left(\frac{\pi}{\theta(e^t)}\right)^2\mu(t),
\end{equation}
where $\theta(r)$ is defined in the following way.
If $u(z)>0$ on the circle $\{ z:|z|=r\}$ then $\theta(r)=+\infty$;
otherwise $\theta(r)$ is the
angular measure of the set $\{ \phi: u(re^{i\phi})>0\}.$

Inequality (\ref{carleman}) expresses a kind of convexity.
If $\theta\equiv+\infty$, then this is the ordinary convexity of $\mu(t)$.
Small $\theta$ implies that $\mu$ is ``very convex''.

It follows that when the set $\{ z\in A:u(z)>0\}$ is narrow,
the function $u$ must grow fast.
When $er_0<r_1\leq+\infty$,
one can derive from (\ref{carleman}) that 
\begin{equation}\label{carl2}
\log m_2(er)\geq\pi\int_{r_0}^r\frac{ds}{s\theta(s)}+\left.\log\frac{dm_2}{d\log r}
\right\vert_{r=r_0},\quad r_0<r<r_1/e.
\end{equation}
It is this version of Carleman's inequality that is most frequently used.
When $r_1=\infty$, inequality (\ref{carl2}) is asymptotically best possible
when $r\to\infty$; extremal regions are angular sectors.

An appropriate version of Carleman's inequality holds in every dimension.
\vspace{.1in}

To derive Theorem 1 from Carleman's inequality we suppose without loss
of generality that all direct singularities lie over finite points,
and that there are at least two of them.
Choose $p\geq 2$ direct
singularities of $f^{-1}$ and consider Euclidean discs $V_j, 1\leq j\leq p$
of radii $\epsilon$ around the corresponding asymptotic values $a_j$.
We choose $\epsilon$ so small that the sets $D_j=S_j(V_j)$ are
disjoint and $f(z)\neq a_j$ in $D_j$. The last condition can be satisfied
by definition of direct singularity.

Then we consider functions
$$u_j(z)=\log\frac{\epsilon}{|f(z)-a_j|},\quad z\in D_j.$$
Each $u_j$ is positive and harmonic in $D_j$ and zero on the boundary.
We denote 
$$B_j(r)=\max_{|z|=r}u_j(z),\quad B(r)=\max_jB_j(r).$$
Then (\ref{carl2}) implies for every $j$
$$\log B(er)\geq\log m_2(er,u_j)\geq
\pi\int_{r_0}^r\frac{ds}{s\theta_j(s)}+O(1)$$
Averaging in $j$ gives
\begin{equation}\label{27}
\log B(er)\geq
\frac{\pi}{p}\int_{r_0}^r\frac{ds}{s}
\sum_{j=1}^p\frac{1}{\theta_j(s)}+O(1).
\end{equation}
On the other hand, by Cauchy--Schwarz inequality,
$$p^2=\left(\sum_j\sqrt{\theta_j(s)}\frac{1}{\sqrt{\theta_j(s)}}\right)^2\leq
\sum_j\theta_j(s)\sum_j\frac{1}{\theta_j(s)}\leq
2\pi\sum_j\frac{1}{\theta_j(s)}.$$
Inserting this to (\ref{27}), we obtain
$$\log B(er)\geq\frac{p}{2}\int_{r_0}^r\frac{ds}{s}+O(1)=
\frac{p}{2}\log r+O(1).$$
By rudimentary Nevanlinna theory, this implies that 
$$\liminf_{r\to\infty}r^{-p/2}T(r,f)>0,$$
so the lower order of $f$,
is at least $p/2$. When $p=1$ we do not obtain any estimate because
many circles $|z|=r$ may lie entirely in $D_1$.  
\vspace{.1in}

\noindent
{\bf Theorem 2.} (\cite{BE,H}) {\em For a meromorphic
function $f$ of finite lower order, each indirect singularity over
a point $a$
is a limit of critical points whose critical values are distinct from $a$.}
\vspace{.1in}

For functions of infinite order, this is not so; there are
entire functions of infinite order without critical points at all,
and such that the set $\Sing(f)$ has the power of continuum \cite{V},
and all but countably many singularities are indirect.

The main analytic tool in the proof of Theorem 2 is
the Carleman
inequality.

Theorem 2 helps to prove in many situations the existence of
critical points.
The simplest example is the following result, originally established
by Clunie, Eremenko, Langley and Rossi:
\vspace{.1in}

\noindent
{\bf Theorem 3.} {\em Let $f$ be a transcendental meromorphic function
of order $\rho$.

a) If $\rho<1$, then $f'$ has infinitely many zeros,

b) If $\rho<1/2$, then $f'/f$ has infinitely many zeros,

c) If $f$ is entire, and $\rho<1$, then $f'/f$ has infinitely many zeros.
}
\vspace{.1in}

Here is another application of Theorem 2: {\em If $f$ is a non-constant
meromorphic function, then $ff'$ takes every finite non-zero value.}

This was conjectured by Hayman in 1967, and all known proofs of
this conjecture use Theorem 2\footnote{This was written in 2013. It is shown
in \cite{An} that Hayman's conjecture also follows from a deep result
of Yamanoi \cite{Ya}.}

\vspace{.2in}

\noindent
{\bf 4. The sets of singularities of various types}
\vspace{.1in}

The set of asymptotic values is an analytic (Suslin) set. This is a larger
class than Borel sets. One of the several equivalent definitions
is that a Suslin set is a continuous image of a Borel set.

More generally, the set of projections of singularities
of an arbitrary multi-valued analytic function is an analytic set.
This was proved by Mazurkiewicz who introduced his metric specially
for this purpose.

Nothing more can be said,
even if one considers asymptotic values of meromorphic
functions of restricted growth.
\vspace{.1in}

\noindent
{\bf Theorem 3.} (Cant\'on--Drasin--Granados) {\em For every analytic (Suslin)
set $A$,
and every $\lambda\geq 0$ there exists a meromorphic function of order
$\lambda$ whose set of asymptotic values is equal to $A$.}
\vspace{.1in}

This is a difficult result which improves on two earlier 
simpler constructions:

{\em For every analytic set, there is an entire function whose set of
asymptotic values is $A\cup\{\infty\}$} (M. Heins),

and

{\em For every $\lambda\geq 0$, there is a meromorphic function of order 
$\lambda$ whose set of asymptotic values is $\bC$.} (A. Eremenko).

Thus there is no restriction on the size of the set of asymptotic values
of a meromorphic function of given order,
and the Corollary from Theorem 1
is the only restriction for entire functions
of given order.

The asymptotic values coming from
direct singularities are rare:
\vspace{.1in}

\noindent
{\bf Theorem 4.} (M. Heins) {\em Let $f$ be a meromorphic function in $\C$,
$V$ a disc, and $D$ a component of $f^{-1}(V)$.
Then the restriction $f:D\to V$ takes every value in $V$ with at most one
exception.
Thus the set of projections of direct singularities is at most countable.}
\vspace{.1in}

On the other hand, the set of direct singularities over one point
can have the power of continuum \cite{BE2}.

If $a$ is an omitted value of $f$, then there is at least one 
singularity over $a$ (Iversen's theorem). Evidently, all singularities
over $a$ are direct in this case, but it is possible that none of them is
logarithmic.

Non-logarithmic direct singularities
of inverses of {\em entire} functions have additional
interesting properties:
\vspace{.1in}

\noindent
{\bf Theorem 5.} (Sixsmith) {\em Let $a\in\bC$ be a projection
of a direct non-logarithmic
singularity of the inverse of an entire function. Then either $a$ is a limit
of critical values, or every neighborhood of this singularity 
contains another transcendental singularity which is either indirect
or logarithmic and whose projection is different from $a$.}
\vspace{.1in}

\noindent
{\bf Theorem 6.} \cite{BE2} {\em Let $S$ be a direct non-logarithmic
singularity of the inverse of an entire function over a point $a\in\C$. Then
every neighborhood of $S$ contains other direct singularities over the
same point $a$.}
\vspace{.1in}

It follows that whenever the inverse of an entire function has a direct
non-logarithmic singularity over a finite point, it must have
the set of singularities of the power of continuum over the same point.

Another corollary is that direct singularities over finite points
of inverses of entire functions of
finite order are all logarithmic.
\vspace{.2in}

\noindent
{\bf 5. Classes of functions defined by restrictions on their singular values}
\vspace{.2in}

By {\em singular values} we mean critical and asymptotic values.

The simplest class is the Speiser class $\S$ which consists of
meromorphic functions with finitely many critical and asymptotic values.
It is the union of classes
$\S_q$ which consist of functions with $q$ critical and asymptotic values.

These are the simplest meromorphic functions
from the geometric point of view. Examples are $\exp(z),\cos z,\wp(z)$,
rational functions. The class of entire functions in $\S$ is closed under
composition. To describe the most important 
property of functions of class $\S$ we need a definition.
\vspace{.1in}

\noindent
{\em Definition.} {\em Two meromorphic functions $f$ and $g$ defined
in simply connected regions are called
{\em topologically equivalent} if there exist homeomorphisms $\phi$ and $\psi$
such that $f\circ\phi=\psi\circ g.$}
\vspace{.1in}

\noindent
{\bf Theorem 7.} (Teichm\"uller, Eremenko--Lyubich)
{\em Let $f$ be a function of class $\S$, and $M_f$
the set of meromorphic functions equivalent to $f$. Then
all functions in $M_f$ are defined in $\C$, and
the set $M_f$
is a complex analytic manifold
of dimension $q+2$, on which the critical and asymptotic values
are holomorphic. }
\vspace{.1in}

These manifolds were introduced and studied in \cite{EL}.
The crucial property here is that $M_f$ has finite dimension.
This permits to extend Sullivan's proof of the absence of wandering
domains to functions of class $\S$.

In the study of holomorphic families of entire functions, the dependence
of periodic points on parameter is important. Let us fix some $g\in \S$,
and let $f\in M_g$. Consider the equation for a periodic point
$$f^m(z)=z.$$
Solution of this equation $z=\alpha(f)$ is a multi-valued function on $M_g$.
The main result on this multi-valued function is 
\vspace{.1in}

\noindent
{\bf Theorem 8.} (Eremenko--Lyubich) {\em The function
$\alpha$ has only algebraic singularities on $M_g$.} 
\vspace{.1in}

The proof uses a version of Carleman's inequality (\ref{carleman}) which
takes into account not only the ``width'' of $D$ but also the amount
of ``spiraling'' of $D$ as $z\to\infty$.
\vspace{.1in}

Class $\B$ consists of transcendental
entire functions whose set of singular values
is bounded. As $\infty$ is always an asymptotic value of an entire function,
it follows that for $f\in \B$, all singularities over $\infty$
are isolated and thus logarithmic. This means that the behavior of the
inverse $f^{-1}$ near infinity is as simple as possible for a transcendental
function. Evidently $\S\subset\B$.

An important analytic tool in the study of functions of
class $\B$ is the following
\vspace{.1in}

\noindent
{\bf Theorem 9.} (Eremenko--Lyubich)
{\em For every $f\in\B$, there exists
$R>0$ such that whenever $|f(z)|>R$,
we have
$$\left|z\frac{f'(z)}{f(z)}
\right|\geq\frac{1}{4\pi}\left(\log|f(z)|-\log R\right).$$
}
\vspace{.1in}

Recently, Sixsmith found that this property actually characterizes 
class~$\B$:
\vspace{.1in}

\noindent
{\bf Theorem 10.} {\em Let $f$ be a transcendental entire function. Then
either $f\in\B$ and
$$\eta:=\lim_{R\to\infty}\inf_{|f(z)|>R}\left|z\frac{f'(z)}{f(z)}\right|
=+\infty,$$
or $f\not\in\B$ and $\eta=0$.}
\vspace{.1in}

The proof of Theorem 10 uses Theorem 5.

Individual functions of class $\S$ were studied from the point
of view of the general theory of meromorphic functions
by Nevanlinna, Teichm\"uller and others. It was found in \cite{EL}
that classes $\S$ and $\B$ are interesting from the point of view of
dynamics.

In general, entire functions, can have all sorts of dynamical
pathology: 
they can have wandering domains, measurable invariant line fields on the
Julia set, and invariant components of the set of normality which do not
contain singular values, and where the iterates converge to infinity
(Baker domains).

The proof of Sullivan's non-wandering theorem depends on the fact that
rational functions topologically conjugate to a given function form a 
manifold of finite dimension. Thus Theorem~{5} permits to extend
this result to the class $\S$, essentially with the same proof.

Other pathologies in dynamics of entire functions
are apparently related to the complicated behavior near $\infty$.
So the class $\B$ with simplest possible behavior
near infinity was introduced
and Theorem~6 can be used to show that
for functions of class $\B$ the iterates cannot
tend to infinity on the set of normality. So Baker domains do not
exist for such functions.

These results permitted to obtain 
a classification of periodic components of the set of normality
for functions of class $\S$, similar to such classification
for rational functions. 

Since then, these classes were intensively studied in holomorphic
dynamics. We mention only one recent result of Bishop: 
{\em functions of class $\B$ can have wandering domains}.
\vspace{.2in}

\noindent
{\bf 6. Further applications to dynamics}
\vspace{.1in}

Let $f$ be a meromorphic function, and $z_0$ a periodic point
which means that $f^nz_0=z_0$ for some $n$.
The smallest $n$ with this property
is called the order of $z_0$, and the derivative $\lambda=(f^n)^\prime(z_0)$,
where $n$ is the order,
is called the multiplier. A periodic point is called attracting, repelling
or neutral if $|\lambda|<1,\; |\lambda|>1$ or $|\lambda|=1$, respectively.
If $\lambda^m=1$ for some $m$, the periodic point is called neutral rational.

If $z_0$ is a periodic point of order $n$ then $z_0$ is a fixed point
of $f^n$ with the same multiplier. Attracting fixed points and neutral fixed
points with multiplier $1$ have non-empty immediate basins of attraction.
An immediate basin of attraction is defined as a maximal invariant
region where $f^m(z)\to z_0$ as $n\to\infty$.
The following theorem of Fatou is fundamental in holomorphic dynamics:
\vspace{.1in}

\noindent
{\bf Theorem 11.} {\em Let $D$ be an immediate basin of attraction
of an attracting or of a neutral point with multiplier $1$.
Then the restriction $f:D\to D$ cannot be a covering map.
This means that this map has either a critical point of an asymptotic
curve with asymptotic value in $D$. Moreover, the trajectory of some
singular value in $D$ is not absorbed by $z_0$.}
\vspace{.1in}

In the simplest case that $z_0$ is attracting, there is a one-line proof
of the first part of the theorem. Suppose that $f:D\to D$ is a covering.
Then $f$ is a local isometry with respect to the hyperbolic metric in $D$.
But at $z_0$, $f$ strictly compresses the hyperbolic metric. This contradiction
proves the statement.

The main corollary of Theorem 11 is that for rational functions,
or for functions of class $\S$, the number of attracting and neutral rational
cycles is finite: it does not exceed the number of singular values.

\vspace{.2in}

A component $D$ of the set of normality
is called completely invariant if $f^{-1}(D)=D$.
The boundary of such component must coincide with the Julia set.
It follows that if a meromorphic function $f$ has at least two
completely invariant components, they all must be simply connected.

If $f$ is a rational function with at least two completely invariant
components $D_j$, then $f:D_j\to D_j$ are ramified coverings
of degree $d=\deg f$,
so by the Riemann--Hurwitz theorem, each $D_j$
must contain $d-1$ critical points
(counting with multiplicity). It follows that a rational function can have
at most two completely invariant components.

How many completely invariant components
can an entire meromorphic function have, is not known\footnote{Baker's proof
\cite{Ba} that a transcendental entire function has at most one completely
invariant domain contains a gap, \cite{RS}.}.
It is conjectured that the answer is at most one for transcendental entire
functions and 
at most two for meromorphic functions, and this is known in the case of meromorphic
functions of class $\S$, see \cite{BKL}.
For meromorphic functions of class $\S$ with
two completely invariant components,
the Julia set is a Jordan curve \cite{BE1}.
\vspace{.2in}

\noindent
{\bf Appendix. Proof of Carleman's inequality}
\vspace{.2in}

Let $v(t,\phi)=u(e^{t+i\phi})$; this is a sub harmonic function,
and we suppose for simplicity that it is continuous.
Let $D=\{ z:u(z)>0\}$.
All integrals below are over the arcs
$\{\phi:e^{t+i\phi}\in D\},$ and we omit $d\phi$.
We have
$$\nu(t):=\int v^2=\mu^2(t),$$
$$\nu^\prime=2\int vv_t,$$
\begin{equation}\label{a}
\nu^{\prime\prime}=2\int(v_t^2+vv_{tt})\geq 2\int(v_t^2+v_\phi^2),
\end{equation}
where we used $u_{tt}+u_{\phi\phi}\geq 0$ and integrated by parts.

Wirtinger's inequality for a $C^1$ function which equals
to zero at the endpoints of an interval $I$
says
$$|I|^2\int_Iv^2\geq\pi^2\int_I v_\phi^2.$$
Applying this to each maximal interval where $v>0$ we obtain
$$\int v_\phi^2\geq\left(\frac{\pi}{\ell}\right)^2\int v^2,$$
where $\ell(t)=\theta(e^t)$.
Cauchy's inequality gives
$$(\nu^\prime)^2=4\left(\int vv_t\right)^2\leq 4\int v^2\int v_t^2=4\nu\int v_t^2.$$
Combining these two inequalities with (\ref{a}) we obtain
$$\nu^{\prime\prime}\geq\frac{(\nu^\prime)^2}{2\nu}+2\left(\frac{\pi}{\ell}\right)^2\nu.$$
Rewriting this for $\mu=\sqrt{\nu}$ we obtain
$$\mu^{\prime\prime}\geq\left(\frac{\pi}{\ell}\right)^2\mu,$$
which is (\ref{carleman}).

To obtain (\ref{carl2}) we set $\omega=\log\mu$. Then
$$\omega^\prime=\frac{\mu^\prime}{\mu},\quad\omega^{\prime\prime}=
\frac{\mu^{\prime\prime}}{\mu}-\left(\frac{\mu^\prime}{\mu}\right)^2,$$
so 
$$\omega^{\prime\prime}+(\omega^\prime)^2=\frac{\mu^{\prime\prime}}{\mu}\geq\left(\frac{\pi}{\ell}\right)^2.$$
Now
$$\left(\omega^\prime+\frac{\omega^{\prime\prime}}{2\omega^\prime}\right)^2\geq
(\omega^\prime)^2+\omega^{\prime\prime}\geq\left(\frac{\pi}{\ell}\right)^2.$$
Thus
$$\omega^\prime+
\frac{\omega^{\prime\prime}}{2\omega^\prime}\geq\frac{\pi}{\ell}.$$
Now notice that 
$$\omega^\prime+\frac{\omega^{\prime\prime}}{2\omega^\prime}=
\frac{1}{2}\frac{d}{dt}\log\left(\frac{d}{dt}e^{2\omega}\right),$$
thus
$$\frac{d}{dt}\log\left(\frac{d\nu}{dt}\right)\geq\frac{2\pi}{\theta(e^t)}.$$
Returning to the variable $r=e^t$ and
integrating this twice with respect to $t$ we obtain (\ref{carl2}).

{\em Purdue University, West Lafayette IN 47907 USA

eremenko@math.purdue.edu}

\end{document}